\documentclass[runningheads]{lncse}

\usepackage{amsmath}
\usepackage{amsfonts}
\usepackage{epsfig}
\usepackage{graphics} 
\usepackage{hyperref}

\usepackage{color}
\usepackage{url}
\usepackage{bm}

\usepackage{xifthen}
\usepackage{stmaryrd}

\newcommand{\primal}{w}
\newcommand{\primalt}{v}

\newcommand{\aux}{\bm{M}}
\newcommand{\auxt}{\bm{N}}

\newcommand{\stiffnessC}{\mathcal{C}}

\newcommand{\vp}{p}
\newcommand{\vpt}{q}
\newcommand{\vphi}{\phi}
\newcommand{\vphit}{\psi}

\newcommand{\lagrange}{\lambda}

\newcommand{\lagranget}{\mu} 
\newcommand{\lagrangetNDer}{\mu_n}
\newcommand{\lagrangetTDer}{\mu_t}

\newcommand{\Q}{{Q}}
\newcommand{\Vprimal}{W}

\newcommand{\VRegDecomp}{V}
\newcommand{\SpPsi}{\Psi}
\newcommand{\SpLagrange}{\Lambda}

\newcommand{\Htwob}[1][]{
\ifthenelse{\isempty{#1}}
{\bm{H}^2(\Omega)}
{\bm{H}^2_0(\Omega)}
}

\newcommand{\FE}{\mathcal{S}}

\newcommand{\GammaVertices}[1][]{
\ifthenelse{\isempty{#1}}
{\mathcal{V}_\Gamma}
{\mathcal{V}_{\Gamma,#1}}
}
\newcommand{\GammaEdges}[1][]{
\ifthenelse{\isempty{#1}}
{\mathcal{E}_\Gamma}
{\mathcal{E}_{\Gamma,#1}}
}

\newcommand{\vertices}[1][]{
\ifthenelse{\isempty{#1}}
{\mathcal{V}_h}
{\mathcal{V}_{h,#1}}
}
\newcommand{\edges}[1][]{
\ifthenelse{\isempty{#1}}
{\mathcal{E}_h}
{\mathcal{E}_{h,#1}}
}

\newcommand{\vertexJump}[1]{\llbracket #1 \rrbracket}

\newcommand{\grad}{\nabla}

\renewcommand{\div}{\operatorname{div}}

\newcommand{\bdiv}{\operatorname{Div}}
\newcommand{\bCurl}{\operatorname{Curl}}

\newcommand{\symbCurl}{\operatorname{symCurl}}

\newcommand{\Id}{\bm{I}}

\newcommand{\matrixTrace}{\mathrm{tr}}

\newcommand{\extPsi}[1]{\vphit[#1]}

\DeclareMathOperator{\dotprod}{\cdot}
\newcommand{\pp}[2]{\partial_{#2} #1}

\begin{document}

\title{On a new mixed formulation of Kirchhoff plates on curvilinear polygonal domains}

\titlerunning{Mixed formulation of Kirchhoff plates}

\author{Katharina Rafetseder \and Walter Zulehner}

\authorrunning{Katharina Rafetseder et al.}   

\institute{
Johannes Kepler University Linz, Institute of Computational Mathematics, Altenberger Stra{\ss}e 69, 4040 Linz, Austria {\tt \{rafetseder,zulehner\}@numa.uni-linz.ac.at}
}

\maketitle

\begin{abstract}
For Kirchhoff plate bending problems on domains whose boundaries are curvilinear polygons a discretization method based on the consecutive solution of three second-order problems is presented.

In Rafetseder and Zulehner (preprint, arXiv:1703.07962) a new mixed variational formulation of this problem is introduced using a nonstandard Sobolev space (and an associated regular decomposition) for the bending moments. In case of a polygonal domain the coupling condition for the two components in the decomposition can be interpreted as standard boundary conditions, which allows for an equivalent reformulation as a system of three (consecutively to solve) second-order elliptic problems. 

The extension of this approach to curvilinear polygonal domains poses severe difficulties. Therefore, we propose in this paper an alternative approach based on Lagrange multipliers. 
\end{abstract}

\section{The Kirchhoff plate bending problem}
We consider the Kirchhoff plate bending problem, where the undeformed mid-surface is described by a domain $\Omega \subset \mathbb{R}^2$ with a Lipschitz boundary $\Gamma$. The plate is considered to be clamped on a part $\Gamma_c \subset \Gamma$, simply supported on $\Gamma_s \subset \Gamma$, and free on $\Gamma_f \subset \Gamma$  with $\Gamma = \Gamma_c \cup \Gamma_s \cup \Gamma_f$. Furthermore, $n = (n_1,n_2)^T$ and $t= (-n_2,n_1)^T$ represent the unit outer normal vector and the unit counterclockwise tangent vector to $\Gamma$, respectively. 

Then the problem reads: For given load $f$, find a deflection $\primal$ such that
\begin{equation}\label{rafetseder_contrib:eq:class_formulation}
	\div\bdiv \big(\stiffnessC  \grad^2 \primal\big) =  f \quad \text{in } \Omega,
\end{equation}
where $\div$ denotes the standard divergence of a vector-valued function, $\bdiv$ the row-wise divergence of a matrix-valued function, $\grad^2$ the Hessian, and $\stiffnessC$ a fourth-order material tensor.
The boundary conditions are given by
\begin{equation*}
\begin{alignedat}{4} 
	&\primal = 0, \quad &&\pp{\primal}{n}=0 & \quad \text{on} \ \Gamma_c, \\
	&\primal = 0, \quad &&\aux n \dotprod n = 0 & \quad \text{on} \ \Gamma_s, \\
	&\aux n \dotprod n = 0, \quad &&\pp{(\aux n \dotprod t)}{t} + \bdiv \aux \dotprod n = 0 & \quad \text{on} \ \Gamma_f,
\end{alignedat}
\end{equation*}
and the corner conditions
\begin{equation*}
	\vertexJump{\aux_{nt}}_x = (\aux n_1 \dotprod t_1)(x) - (\aux n_2 \dotprod t_2)(x) = 0  \quad \text{for all} \ x \in \GammaVertices[f],
\end{equation*}
where $\aux$ denotes the bending moment tensor, given by $\aux = -\stiffnessC \grad^2 \primal$, and $\GammaVertices[f]$ denotes the set of corner points whose two adjacent edges (with corresponding normal and tangent vectors $n_1$, $t_1$ and $n_2$, $t_2$) belong to $\Gamma_f$.

As an example, the material tensor $\stiffnessC$ for isotropic materials is given by
\begin{equation}\label{rafetseder_contrib:eq:materialTensor}
	\stiffnessC \auxt = D \bigl( (1-\nu)\auxt + \nu \ \matrixTrace(\auxt) \Id \bigr),
\end{equation}
for matrices $\auxt$, where $\nu$ is the Poisson ration, $D>0$ depends on $\nu$, Young's modulus, and the thickness of the plate, $\Id$ is the identity matrix, and $\matrixTrace$ is the trace operator for matrices.

A standard (primal) variational formulation of \eqref{rafetseder_contrib:eq:class_formulation} is given as follows: Find $\primal \in \Vprimal$ such that
\begin{equation}\label{rafetseder_contrib:eq:primal_formulation}
	\int_\Omega \stiffnessC \grad^2 \primal : \grad^2 \primalt \ dx = \langle F,\primalt \rangle \quad \text{for all} \ \primalt \in \Vprimal,
\end{equation}
with the Frobenius inner product $\bm{A} : \bm{B} = \sum_{i,j} \bm{A}_{ij} \, \bm{B}_{ij}$ for matrices $\bm{A}, \bm{B}$, the right-hand side $\langle F,\primalt \rangle = \int_\Omega f \, v \ dx$, and the function space
\begin{equation}
	\Vprimal = \{\primalt \in H^2(\Omega) : \ \primalt = 0, \ \pp{\primalt}{n} = 0 \ \text{on} \ \Gamma_c, \quad \primalt = 0 \ \text{on} \ \Gamma_s \}.
\end{equation}
Here and throughout the paper  $L^2(\Omega)$ and $H^m(\Omega)$ denote the standard Lebesgue and Sobolev spaces of functions on $\Omega$ with corresponding norms $\|.\|_{0}$ and $\|.\|_{m}$ for positive integers $m$. Moreover, $H^1_{0,\Gamma'}(\Omega)$ denotes the set of function in $H^1(\Omega)$ which vanish on a part $\Gamma'$ of $\Gamma$.
The $L^2$-inner product on $\Omega$ and $\Gamma'$ are always denoted by $(.,.)$ and $(.,.)_{\Gamma'}$, respectively, no matter whether it is used for scalar, vector-valued, or matrix-valued functions. We use $H^*$ to denote the dual of a Hilbert space $H$ and $\langle .,. \rangle$ for the duality product on $H^*\times H$.

\section{New mixed formulation}
In our previous work \cite{rafetseder_zulehner_2017} a new mixed variational formulation for the Kirchhoff plate bending problem with the bending moment tensor $\aux$ as additional unknown is derived. The new mixed formulation satisfies Brezzi's conditions and is equivalent to the original problem without additional convexity assumption on $\Omega$. These important properties come at the expense of an appropriate nonstandard Sobolev space for $\aux$.
In order to make this space computationally accessible, we show in \cite[Theorem 4.2]{rafetseder_zulehner_2017} a regular decomposition of it, which provides the following representation of the solution $\aux$
\begin{equation*}
	\aux= \vp \Id + \symbCurl \vphi,
\end{equation*}
with $\vp\in \Q = H^1_{0,\Gamma_c \cup \Gamma_s}(\Omega)$ and $\vphi \in (H^1(\Omega))^2$ satisfying the coupling condition
\begin{equation}\label{rafetseder_contrib:eq:coupling_cond_p}
	\langle \pp{\vphi}{t} ,\grad \primalt \rangle_\Gamma = - \int_\Gamma \vp \ \pp{\primalt}{n} \ ds \quad \text{for all} \ \primalt \in \Vprimal,
\end{equation}
where $\partial_t \vphi = (\bCurl \vphi) n \in (H^{-\frac{1}{2}}(\Gamma))^2$ with $H^{-\frac{1}{2}}(\Gamma) = (H^{\frac{1}{2}}(\Gamma))^*$.
Here the symmetric $\bCurl$ is defined as $\symbCurl \vphit = \frac{1}{2}( \bCurl \vphit + (\bCurl \vphit)^T)$ with
\begin{equation*}
\bCurl \vphit=
	\begin{pmatrix}
	\pp{\vphit_{1}}{2} & -\pp{\vphit_{1}}{1}  \\
	\pp{\vphit_{2}}{2} & -\pp{\vphit_{2}}{1}  \\
	\end{pmatrix}.
\end{equation*}

The analogous representation for the test functions associated to $\aux$ in the mixed formulation leads to the following equivalent formulation of \eqref{rafetseder_contrib:eq:primal_formulation}:
For $F\in \Q^*$, find $(\vp, \vphi) \in \VRegDecomp$ and $\primal \in \Q$ such that
\begin{equation}\label{rafetseder_contrib:eq:coupling_cond_formulation}
\begin{alignedat}{3} 
	 & (\vp \Id + \symbCurl \vphi, \vpt \Id + \symbCurl \vphit)_{\stiffnessC^{-1}} & & - (\grad \primal, \grad \vpt) & & = 0,\\
	 & - (\grad \vp, \grad \primalt) & & & & = -\langle F, \primalt \rangle,
\end{alignedat}
\end{equation}
for all $\primalt \in \Q = H^1_{0,\Gamma_c \cup \Gamma_s}(\Omega)$ and $(\vpt, \vphit) \in \VRegDecomp$, where the function space $\VRegDecomp$ is given as the subset of $ (\vpt, \vphit) \in Q \times (H^1(\Omega))^2$ satisfying
\begin{equation}\label{rafetseder_contrib:eq:coupling_cond_q}
	\langle \pp{\vphit}{t} ,\grad \primalt \rangle_\Gamma = - \int_\Gamma \vpt \ \pp{\primalt}{n} \ ds \quad \text{for all} \ \primalt \in \Vprimal.
\end{equation}
Here, we use the notation $(\aux,\auxt)_{\stiffnessC^{-1}} = (\stiffnessC^{-1} \aux,\auxt)$.

\subsection{Coupling condition as standard boundary conditions for $\vphi$}\label{rafetseder_contrib:subsec:boundary_conditions}
In \cite{rafetseder_zulehner_2017} we obtain for polygonal domains $\Omega$ an equivalent formulation of the Kirchhoff plate bending problem \eqref{rafetseder_contrib:eq:primal_formulation} in terms of three (consecutively to solve) second-order elliptic problems:
\begin{enumerate}
\item
The $\vp$-problem: Find $p \in \Q$ such that
\begin{equation*}
	(\grad \vp, \grad \primalt) = \langle F, \primalt \rangle \quad \text{for all} \ \primalt \in \Q .
\end{equation*}
\item
The $\vphi$-problem: For given $\vp \in \Q$, find  $\vphi \in \SpPsi_\vp = \extPsi{\vp} + \SpPsi_0$ such that
\begin{equation*}
   (\symbCurl \vphi, \symbCurl \vphit_0)_{\stiffnessC^{-1}}  = - (\vp \Id, \symbCurl \vphit_0)_{\stiffnessC^{-1}} \quad \text{for all} \ \vphit_0 \in \SpPsi_0 .
\end{equation*}
\item
The $\primal$-problem: For given $\aux = \vp \Id + \symbCurl \vphi$, find $\primal \in \Q$ such that
\begin{equation*}
	(\grad \primal, \grad \vpt) = (\aux, \vpt \Id + \symbCurl \extPsi{\vpt} )_{\stiffnessC^{-1}} \quad \text{for all} \ \vpt \in \Q.
\end{equation*}
\end{enumerate}
The second and the third problem require the construction of a particular function $\extPsi{\vpt}$ satisfying the coupling condition \eqref{rafetseder_contrib:eq:coupling_cond_q} for given $\vpt\in\Q$, for details see \cite{rafetseder_zulehner_2017}. The space $\SpPsi_0$ consists of all functions in $(H^1(\Omega))^2$ satisfying \eqref{rafetseder_contrib:eq:coupling_cond_q} for $\vpt = 0$.

The approach presented in \cite{rafetseder_zulehner_2017} is to characterize $\SpPsi_0$ as space of functions $\vphit$ with standard boundary conditions available in $(H^1(\Omega))^2$. 
Originally, the boundary conditions for  $\vphit \in \SpPsi_0$ are, roughly speaking, conditions for tangential derivatives of  $\vphit$ of the form
\begin{align}
	&\pp{\vphit}{t} \dotprod n = 0 \quad && \text{on} \ \Gamma_s,\label{rafetseder_contrib:eq:boundaryCond_phi_s}\\
	&\pp{^2\vphit}{t} \dotprod t  = 0,  \quad \pp{\vphit}{t} \dotprod n = 0 \quad && \text{on} \ \Gamma_f.\label{rafetseder_contrib:eq:boundaryCond_phi_f}
\end{align}
For polygonal domains we obtain from \eqref{rafetseder_contrib:eq:boundaryCond_phi_f} a Dirichlet boundary condition for $\vphit$. Moreover, \eqref{rafetseder_contrib:eq:boundaryCond_phi_s} yields a Dirichlet boundary condition for the normal component $\vphit \dotprod n$.
However, the considerations heavily rely on a polygonal domain and it is not clear how to obtain standard boundary conditions in the curved case. This is our main motivation to investigate an alternative approach to incorporate the coupling condition \eqref{rafetseder_contrib:eq:coupling_cond_q} based on Lagrange multipliers, which we introduce in the next section.

In \cite{rafetseder_zulehner_2017} we propose a discretization method for the above introduced formulation using a Nitsche method to incorporate the boundary conditions in the $\vphi$-problem and present a numerical analysis of the method.

\section{Coupling condition via Lagrange multipliers}
We consider a domain $\Omega$, whose boundary is a curvilinear polygon of class $C^\infty$. This means that $\Gamma = \bigcup_{k=1}^K \overline E_k$, where the edges $E_k$ are $C^\infty$ curves for $k=1,2, \dots, K$ and $\overline{E}_k$ denotes the closure of $E_k$. The edges are numbered consecutively in counterclockwise direction. We denote the vertex at the endpoint of $\overline E_k$ by $a_k$ and the interior angle at $a_k$ by $\omega_k$. Note, since we consider a closed boundary curve, the index $k=0$ is in the following always identified with $k=K$.

Furthermore, we assume that each edge $E_k$ is contained in exactly one of the sets $\Gamma_c$, $\Gamma_s$, $\Gamma_f$, and the edges are maximal in the sense that two edges with the same boundary condition do not meet at an angle of $\pi$.

By using the representation $\grad \primalt = (\pp{\primalt}{n}) \, n + (\pp{\primalt}{t}) \, t$ and incorporating the boundary conditions for $\primalt \in \Vprimal$, the coupling condition \eqref{rafetseder_contrib:eq:coupling_cond_p} reads
\begin{equation*}
	(\pp{\vphi}{t}\dotprod n + \vp, \pp{\primalt}{n})_{\Gamma_s \cup \Gamma_f} + (\pp{\vphi}{t}\dotprod t, \pp{\primalt}{t})_{\Gamma_f} = 0 \quad \text{for all} \ \primalt \in \Vprimal,
\end{equation*}
provided $\pp{\phi}{t} \in L^2(\Gamma)$. We can rewrite the condition as follows
\begin{equation}\label{rafetseder_contrib:eq:constraint}
	 \sum_{E_k \subset \Gamma_s \cup \Gamma_f} (\pp{\vphi}{t}\dotprod n \ + \ \vp,  \lagrangetNDer^k)_{E_k} +  \sum_{E_k \subset \Gamma_f} (\pp{\vphi}{t}\dotprod t, \lagrangetTDer^k)_{E_k} = 0,
\end{equation}
for all $\lagranget = ((\lagrangetTDer^1, \lagrangetTDer^2, \dots, \lagrangetTDer^K), (\lagrangetNDer^1, \lagrangetNDer^2, \dots, \lagrangetNDer^K)) \in \SpLagrange$ where 
\begin{equation*}
	\SpLagrange = \{(\pp{\primalt}{t}, \pp{\primalt}{n}) \colon \text{for} \ \primalt \in \Vprimal\},
\end{equation*}
with
\begin{equation*}
	\pp{\primalt}{t} = (\pp{\primalt}{t}|_{E_1}, \pp{\primalt}{t}|_{E_2}, \dots, \pp{\primalt}{t}|_{E_K}), \quad \pp{\primalt}{n} = (\pp{\primalt}{n}|_{E_1}, \pp{\primalt}{n}|_{E_2}, \dots, \pp{\primalt}{n}|_{E_K}).
\end{equation*}

We view the original formulation \eqref{rafetseder_contrib:eq:coupling_cond_formulation} as optimality system with constraint $(\grad \vp, \grad \primalt) = \langle F, \primalt \rangle$ and replace the space $\VRegDecomp$ by $Q \times (H^1(\Omega))^2$ and add \eqref{rafetseder_contrib:eq:constraint} as additional constraint. The corresponding optimality system is the starting point for the discretization method we introduce in Sect.~\ref{rafetseder_contrib:sec:discretization}.

\subsection{Characterization of $\SpLagrange$}
\label{rafetseder_contrib:subsec:characterization}
In this subsection we provide an explicit characterization of $\SpLagrange$. Let us consider $\lagranget = ((\lagrangetTDer^1, \lagrangetTDer^2, \dots, \lagrangetTDer^K), (\lagrangetNDer^1, \lagrangetNDer^2, \dots, \lagrangetNDer^K))$, where $\lagrangetTDer^k$ and $\lagrangetNDer^k$ for $k=1,2, \dots, K$ are Lipschitz continuous functions on $\overline E_k$. Then $\lagranget \in \SpLagrange$ if and only if the following three conditions are satisfied:
\begin{enumerate}
	\item The boundary conditions $\lagrangetTDer^k = 0$ on edges $E_k \subset \Gamma_s \cup \Gamma_c$ and $\lagrangetNDer^k = 0$ on edges $E_k \subset \Gamma_c$ have to hold.
	
	\item On each connected component $C$ of $\Gamma_f$ the compatibility condition 
	\begin{equation*}
		\sum_{E_k \subset C} \int_{E_k} \lagrangetTDer^k \ ds = 0 
	\end{equation*}
	has to be satisfied.
	
	\item The four corner values $\lagrangetTDer^{k-1}(a_k), \lagrangetTDer^{k}(a_k), \lagrangetNDer^{k-1}(a_k), \lagrangetNDer^{k}(a_k)$ have to be coupled appropriately.
	For the case $\omega_k \neq \pi$, the conditions are given by
	\begin{equation}\label{rafetseder_contrib:eq:conditionLagrange}
	\begin{aligned}
		\lagrangetTDer^{k-1} (a_k) + \cos\omega_k \ \lagrangetTDer^k (a_k) - \sin\omega_k \ \lagrangetNDer^k (a_k) &= 0,\\
		\cos\omega_k \ \lagrangetTDer^{k-1} (a_k) + \sin\omega_k \ \lagrangetNDer^{k-1} (a_k) + \lagrangetTDer^k (a_k) &= 0,
	\end{aligned}
	\end{equation}
	for all $k=1,2,\dots,K$. These conditions follow as special case from \cite[Theorem 1.5.2.8]{grisvard_1985}.
	\begin{remark}
		In order to describe a change of boundary condition we may also consider an interior angle $\omega_k$ of $\pi$. A corresponding adaption of the conditions \eqref{rafetseder_contrib:eq:conditionLagrange} can be found in \cite{grisvard_1985}.
	\end{remark}
\end{enumerate}
In the following we fix a corner $a_k$ and work out the relation implied by the corresponding boundary conditions and the conditions \eqref{rafetseder_contrib:eq:conditionLagrange} for the four involved quantities $\lagrangetTDer^{k-1}(a_k)$, $\lagrangetTDer^{k}(a_k)$, $\lagrangetNDer^{k-1}(a_k)$, $\lagrangetNDer^{k}(a_k)$, where we skip in the following the argument $a_k$ for better readability. We distinguish three situations:
\begin{enumerate}
	\item Let $a_k$ be an interior corner point of $\Gamma_f$. Then the conditions \eqref{rafetseder_contrib:eq:conditionLagrange} lead to
	\begin{equation*}
		\lagrangetNDer^{k-1} = -\frac{1}{\sin w_k} (\cos w_k \ \lagrangetTDer^{k-1} +  \lagrangetTDer^k), \quad \lagrangetNDer^{k} = \frac{1}{\sin w_k} (\lagrangetTDer^{k-1} + \cos w_k \ \lagrangetTDer^k),
	\end{equation*}
	for arbitrary $\lagrangetTDer^{k-1}$ and $\lagrangetTDer^{k}$.
	
	\item Let $a_k$ be a corner point on the interface of $E_{k-1} \subset \Gamma_s$ and $E_k \subset \Gamma_f$. Then the conditions \eqref{rafetseder_contrib:eq:conditionLagrange} provide
	\begin{equation*}
		\lagrangetTDer^{k-1} = 0, \quad  \lagrangetNDer^{k-1} = -\frac{1}{\sin w_k} \lagrangetTDer^k, \quad \lagrangetNDer^{k} = \frac{1}{\sin w_k} \cos w_k \ \lagrangetTDer^k,
	\end{equation*}
	where $\lagrangetTDer^{k}$ can be freely chosen.
	For the reverse case $E_{k-1} \subset \Gamma_f$ and $E_k \subset \Gamma_s$ an analogous result holds.
	
	\item In all other cases, we obtain $\lagrangetTDer^{k-1} = \lagrangetTDer^k = \lagrangetNDer^{k-1} = \lagrangetNDer^k = 0$.
\end{enumerate}

\subsection{The discretization method}
\label{rafetseder_contrib:sec:discretization}


Let $\FE_h(\Omega)$ be a finite dimensional subspace of $H^1(\Omega)$ of piecewise polynomials (with respect to a subdivision of $\Omega$) and we set $\FE_{h,0}(\Omega) = \FE_h(\Omega) \cap H^1_{0,\Gamma_c \cup \Gamma_s}(\Omega)$. The restriction of functions from $\FE_h(\Omega)$ to $E_k$ is defined as
\begin{equation*}
	\FE_h(E_k) = \{\primalt|_{E_k} \colon \primalt \in \FE_h(\Omega) \}.
\end{equation*}
The discrete space $\SpLagrange_h$ consists of all $\lagranget_h = ((\lagrangetTDer^1, \lagrangetTDer^2, \dots, \lagrangetTDer^K), (\lagrangetNDer^1, \lagrangetNDer^2, \dots, \lagrangetNDer^K))$, where $\lagrangetTDer^k \in \FE_h(E_k)$ and $\lagrangetNDer^k \in \FE_h(E_k)$ for $k=1,2, \dots, K$, subject to the constraints derived in Sect.~\ref{rafetseder_contrib:subsec:characterization}.

In the discrete setting the original formulation \eqref{rafetseder_contrib:eq:coupling_cond_formulation} is equivalent to three (consecutively to solve) second-order problems:
\begin{enumerate}
	\item The discrete $\vp$-problem: Find $\vp_h \in \FE_{h,0}(\Omega)$ such that
	\begin{equation*}
		(\grad \vp_h, \grad \primalt_h) = \langle F, \primalt_h \rangle \quad \text{for all} \ \primalt_h \in \FE_{h,0}(\Omega).
	\end{equation*}
	\item The discrete ($\vphi$, $\lagrange$)-problem:\\ 
	For given $\vp_h \in \FE_{h,0}(\Omega)$, find $\vphi_h \in (\FE_h(\Omega))^2/RT_0$ and $\lagrange_h \in \SpLagrange_h$ such that
	\begin{equation*}
	\begin{alignedat}{4}
		& (\symbCurl \vphi_h, \symbCurl \vphit_h)_{\stiffnessC^{-1}} & & + l_\vphi(\vphit_h,\lagrange_h) & & = - (\vp_h \Id, \symbCurl \vphit_h)_{\stiffnessC^{-1}}\\
		& l_\vphi(\vphi_h,\lagranget_h) & & & & = -l_\vp(\vp_h, \lagranget_h),
	\end{alignedat}
	\end{equation*}
	for all $\vphit_h \in (\FE_h(\Omega))^2/RT_0$ and $\lagranget_h \in \SpLagrange_h$, where
	\begin{align*}
		l_\vphi(\vphi,\lagranget) &= \sum_{E_k \subset \Gamma_s \cup \Gamma_f} (\pp{\vphi}{t}\dotprod n,  \lagrangetNDer^k)_{E_k} +  \sum_{E_k \subset \Gamma_f} (\pp{\vphi}{t}\dotprod t, \lagrangetTDer^k)_{E_k},\\
		l_\vp(\vp, \lagranget) &= \sum_{E_k \subset \Gamma_f} (\vp,  \lagrangetNDer^k)_{E_k},
	\end{align*}
	for $\lagranget = ((\lagrangetTDer^1, \lagrangetTDer^2, \dots, \lagrangetTDer^K), (\lagrangetNDer^1, \lagrangetNDer^2, \dots, \lagrangetNDer^K))$.
	Here, we use the notation $RT_0 = \{ax + b \colon a \in \mathbb{R}, \ b \in \mathbb{R}^2 \}$.
	\item The discrete $\primal$-problem: For given $\aux_h = \vp_h \Id + \symbCurl \vphi_h$ and $\lagrange_h \in \SpLagrange_h$, find $\primal_h \in \FE_{h,0}(\Omega)$ such that
	\begin{equation*}
		(\grad \primal_h, \grad \vpt_h) = (\aux_h, \vpt_h \Id)_{\stiffnessC^{-1}} + l_\vp(\vpt_h, \lagrange_h) \quad \text{for all} \ \vpt_h \in \FE_{h,0}(\Omega).
	\end{equation*}
\end{enumerate}
In comparison with the decoupled formulation in Sect.~\ref{rafetseder_contrib:subsec:boundary_conditions}, here the second problem, the ($\vphi$, $\lagrange$)-problem, is a saddle point problem.

\section{Numerical tests}

As discretization space $\FE_h(\Omega)$ we consider B-splines of degree $p \geq 1$ with maximum smoothness; see, e.g, \cite{cotrell_hughes_brazilevs_2009,daVeiga_buffa_sangalli_vazquez_2014} for more information on this space in the context of isogeometric analysis (IGA). A sparse direct solver is used for each of the three sub-problems. The implementation is done in the framework of the object-oriented C++ library G+Smo ("Geometry + Simulation Modules") \footnote{\url{https://ricamsvn.ricam.oeaw.ac.at/trac/gismo/wiki/WikiStart}}.

\subsection{Square plate with clamped, simply supported and free boundary}
We consider a square plate $\Omega = (-1, 1)^2$ with simply supported north and south boundaries, clamped west boundary and free east boundary. The material tensor $\stiffnessC$ is given as in \eqref{rafetseder_contrib:eq:materialTensor} with $D=1$, $\nu=0$ and the load is $f(x,y) = 4 \pi ^4 \sin (\pi  x) \sin (\pi  y)$.
The exact solution is written in the form
\begin{equation*}
	\primal(x,y) = \bigl((a + b x)\cosh(\pi x) + (c + d x)\sinh(\pi x) + \sin(\pi x)\bigr) \sin(\pi y),
\end{equation*}
which automatically satisfies the boundary conditions on the simply supported boundary parts. The constants $a, b, c$ and $d$ are chosen such that the four remaining boundary conditions (on the clamped and free boundary parts) are satisfied, for details, see \cite{reddy_2007}.
In Table~\ref{rafetseder_contrib:tab:square_k_1} and Table~\ref{rafetseder_contrib:tab:rectangle_k_3} the discretization errors for $p=1,3$ are presented. The first column shows the refinement level $L$, the next three pairs of columns show the respective discretization error and the error reduction relative the previous level. The results show optimal convergence rates for $\primal$ and $\aux$.

\begin{table}
\caption{Discretization errors for square plate, $p=1$}
\label{rafetseder_contrib:tab:square_k_1}
\begin{tabular}{p{0.5cm}p{1.7cm}p{1.3cm}p{1.7cm}p{1.3cm}p{1.9cm}p{0.8cm}}
\hline\noalign{\smallskip}
$L$ & $\|\primal-\primal_h\|_0$ & order & $\|\primal-\primal_h\|_1$ & order & $\|\aux-\aux_{h}\|_0$ & order\\
4 &  $2.82 \cdot 10^{-2}$ & $1.909$ & $6.83 \cdot 10^{-1}$ & $0.992$ & $2.83 \cdot 10^{0}$  & $0.975$ \\ 
5 &  $7.17 \cdot 10^{-3}$ & $1.976$ & $3.42 \cdot 10^{-1}$ & $0.998$ & $1.42 \cdot 10^{0}$  & $0.993$ \\
6 &  $1.80 \cdot 10^{-3}$ & $1.994$ & $1.71 \cdot 10^{-1}$ & $0.999$ & $7.13 \cdot 10^{-1}$ & $0.998$ \\
7 &  $4.50 \cdot 10^{-4}$ & $1.998$ & $8.55 \cdot 10^{-2}$ & $0.999$ & $3.56 \cdot 10^{-1}$ & $0.999$ \\
\noalign{\smallskip}\hline\noalign{\smallskip}
\end{tabular}
\end{table}
\begin{table}
\caption{Discretization errors for square plate, $p=3$}
\label{rafetseder_contrib:tab:rectangle_k_3}
\begin{tabular}{p{0.5cm}p{1.7cm}p{1.3cm}p{1.7cm}p{1.3cm}p{1.9cm}p{0.8cm}}
\hline\noalign{\smallskip}
$L$ & $\|\primal-\primal_h\|_0$ & order & $\|\primal-\primal_h\|_1$ & order & $\|\aux-\aux_{h}\|_0$ & order\\
4 &  $5.47 \cdot 10^{-5}$ & $4.147$ & $2.75 \cdot 10^{-3}$ & $3.070$ & $1.10 \cdot 10^{-2}$  & $3.104$ \\ 
5 &  $3.41 \cdot 10^{-6}$ & $4.001$ & $3.46 \cdot 10^{-4}$ & $2.993$ & $1.38 \cdot 10^{-3}$  & $2.989$ \\
6 &  $2.15 \cdot 10^{-7}$ & $3.984$ & $4.37 \cdot 10^{-5}$ & $2.985$ & $1.75 \cdot 10^{-4}$ & $2.978$ \\
7 &  $1.35 \cdot 10^{-8}$ & $3.989$ & $5.50 \cdot 10^{-6}$ & $2.989$ & $2.22 \cdot 10^{-5}$ & $2.985$ \\
\noalign{\smallskip}\hline\noalign{\smallskip}
\end{tabular}
\end{table}

\subsection{Circular plate with simply supported boundary}
As a second example, we consider the simply supported circular plate with radius $r=1$ and uniform loading $f=1$. The material tensor $\stiffnessC$ is given as in \eqref{rafetseder_contrib:eq:materialTensor} with $D=1$ and $\nu = 0.3$. The exact solution is given by $\primal(x) = c_1 + c_2 r^2 + c_3 r^4$ where $r^2= x_1^2 + x_2^2$, $c_3=1/64$ and $c_1$, $c_2$ are determined from the boundary conditions. 
For this test reproducing the exact geometry is essential, see the discussion of the so-called Babu\v{s}ka paradox in \cite{babuska_pitkaranta_1990}. Therefore, we use an exact geometry representation by means of non-uniform rational B-splines (NURBS).
In Table~\ref{rafetseder_contrib:tab:disc_k_1} and Table~\ref{rafetseder_contrib:tab:disc_k_3} the discretization errors for $p=1,3$ are presented. The results show optimal convergence rates for $\primal$ and $\aux$.

\begin{table}
\caption{Discretization errors for circular plate, $p=1$}
\label{rafetseder_contrib:tab:disc_k_1}
\begin{tabular}{p{0.5cm}p{1.7cm}p{1.3cm}p{1.7cm}p{1.3cm}p{1.9cm}p{0.8cm}}
\hline\noalign{\smallskip}
$L$ & $\|\primal-\primal_h\|_0$ & order & $\|\primal-\primal_h\|_1$ & order & $\|\aux-\aux_{h}\|_0$ & order\\
4 &  $3.58 \cdot 10^{-4}$ & $1.984$ & $8.37 \cdot 10^{-3}$ & $1.002$ & $8.79 \cdot 10^{-3}$ & $1.020$ \\ 
5 &  $8.98 \cdot 10^{-5}$ & $1.996$ & $4.18 \cdot 10^{-3}$ & $1.000$ & $4.38 \cdot 10^{-3}$ & $1.005$ \\
6 &  $2.24 \cdot 10^{-5}$ & $1.999$ & $2.09 \cdot 10^{-3}$ & $1.000$ & $2.18 \cdot 10^{-3}$ & $1.001$ \\
7 &  $5.62 \cdot 10^{-6}$ & $1.999$ & $1.04 \cdot 10^{-3}$ & $1.000$ & $1.09 \cdot 10^{-3}$ & $1.000$ \\
\noalign{\smallskip}\hline\noalign{\smallskip}
\end{tabular}
\end{table}
\begin{table}
\caption{Discretization errors for circular plate, $p=3$}
\label{rafetseder_contrib:tab:disc_k_3}
\begin{tabular}{p{0.5cm}p{1.7cm}p{1.3cm}p{1.7cm}p{1.3cm}p{1.9cm}p{0.8cm}}
\hline\noalign{\smallskip}
$L$ & $\|\primal-\primal_h\|_0$ & order & $\|\primal-\primal_h\|_1$ & order & $\|\aux-\aux_{h}\|_0$ & order\\
4 &  $4.05 \cdot 10^{-7}$  & $4.319$ & $1.93 \cdot 10^{-5}$ & $3.163$ & $2.01 \cdot 10^{-5}$ & $3.206$ \\ 
5 &  $2.38 \cdot 10^{-8}$  & $4.083$ & $2.35 \cdot 10^{-6}$ & $3.034$ & $2.42 \cdot 10^{-6}$ & $3.054$ \\
6 &  $1.47 \cdot 10^{-9}$  & $4.019$ & $2.93 \cdot 10^{-7}$ & $3.004$ & $3.00 \cdot 10^{-7}$ & $3.013$ \\
7 &  $9.17 \cdot 10^{-11}$ & $4.004$ & $3.67 \cdot 10^{-8}$ & $2.999$ & $3.74 \cdot 10^{-8}$ & $3.003$ \\
\noalign{\smallskip}\hline\noalign{\smallskip}
\end{tabular}
\end{table}

\def\cprime{$'$}

%
%
%
 \bibliographystyle{vmams}  
 \bibliography{KLplate_bibliography}        

\end{document}